\def\ver{Sept. 22, 2003, v.11}
\documentstyle{amsppt}
\magnification=1200
\hsize=6.5truein
\vsize=8.9truein
\hoffset=-10pt
\voffset=-20pt
\topmatter
\title Monodromy Filtration and Positivity
\endtitle
\author Morihiko Saito
\endauthor
\affil RIMS Kyoto University, Kyoto 606-8502 Japan \endaffil
\abstract
We study Deligne's conjecture on the monodromy weight filtration
on the nearby cycles in the mixed characteristic case, and reduce
it to the nondegeneracy of certain pairings in the semistable case.
We also prove a related conjecture of Rapoport and Zink which
uses only the image of the Cech restriction morphism, if
Deligne's conjecture holds for a general hyperplane section.
In general we show that Deligne's conjecture is true if the
standard conjectures hold.
\endabstract
\endtopmatter
\tolerance=1000
\baselineskip=12pt
\def\scirc{\raise.2ex\hbox{${\scriptstyle\circ}$}}
\def\ssbull{\raise.2ex\hbox{${\scriptscriptstyle\bullet}$}}
\def\mopls{\hbox{$\bigoplus$}}
\def\mtim{\hbox{$\times$}}
\def\bC{{\Bbb C}}
\def\bP{{\Bbb P}}
\def\bQ{{\Bbb Q}}
\def\bR{{\Bbb R}}
\def\bZ{{\Bbb Z}}

\def\oE{\overline{E}}

\def\ok{\overline{k}}
\def\oK{\overline{K}}

\def\obQ{\overline{\Bbb Q}}
\def\oal{\overline{\alpha}}
\def\td{\widetilde{d}}
\def\tGam{\widetilde{\Gamma}}
\def\char{\text{\rm char}\,}
\def\Spec{\text{\rm Spec}\,}
\def\NS{\hbox{{\rm NS}}}
\def\CH{\hbox{{\rm CH}}}
\def\Im{\hbox{{\rm Im}}}
\def\Ker{\hbox{{\rm Ker}}}
\def\Hom{\hbox{{\rm Hom}}}
\def\End{\hbox{{\rm End}}}
\def\Tr{\hbox{{\rm Tr}}}
\def\Gal{\hbox{{\rm Gal}}}
\def\Gr{\text{{\rm Gr}}}
\def\Alb{\text{{\rm Alb}}}
\def\prim{\text{\rm prim}}
\def\simto{\buildrel\sim\over\longrightarrow}
\def\SameAuthor{\vrule height3pt depth-2.5pt width1cm}

\document\noindent
\centerline{\bf Introduction}
\footnote""{{\it Date}\,: \ver}

\bigskip\noindent
Let
$ f : X \to S $ be a projective morphism of complex manifolds with
relative dimension
$ n $ where
$ S $ is an open disk and
$ f $ is semistable (i.e.
$ X_{0} := f^{-1}(0) $ is a reduced divisor with normal crossings
whose irreducible components are smooth).
J. Steenbrink [29] constructed a limit mixed Hodge structure by
using a resolution of the nearby cycle sheaf
on which the monodromy weight filtration can be defined.
This limit mixed Hodge structure coincides with the one obtained by
W. Schmid [26] using the
$ \hbox{SL}_{2} $-orbit theorem, because the weight filtration
coincides with the (shifted) monodromy filtration,
see [23] and also [5], [9], [11], [24, 2.3], etc.

A similar construction was then given by M. Rapoport and T. Zink [21]
in the case
$ X $ is projective and semistable over a henselian discrete valuation
ring
$ R $ of mixed characteristic.
Here we may assume
$ f $ semistable because the non semistable case can be reduced to
the semistable case by [2] replacing the discrete valuation field
with a finite extension if necessary.
Then Deligne's conjecture on the monodromy filtration [3, I]
is stated as

\medskip\noindent
{\bf 0.1.~Conjecture.} The obtained weight filtration coincides with
the monodromy filtration shifted by the degree of cohomology.

\medskip
This conjecture is proved so far in the case
$ n \le 2 $ by [21] (using [2]) and in some other cases
(see e.g. [8], [25]).
We cannot apply the same argument as in the complex analytic case,
because we do not have a good notion of positivity for
$ l $-adic sheaves.
In the positive characteristic case, however, the conjecture was
proved by Deligne [4] (assuming
$ f $ is the base change of a morphism to a smooth curve over a
finite field).
In this paper we apply the arguments in [23], [24] to the mixed
characteristic case, and show that Conjecture (0.1) is closely
related to the standard conjectures which would give a notion of
positivity for the pairings of correspondences (but not for the
pairings of cohomology groups in general).

Let
$ Y^{(i)} $ denote the disjoint union of the intersections of
$ i $ irreducible components of the special fiber
$ X_{0} $.
Let
$ k $ be the residue field of
$ R $ which is assumed to be a finite field.
Let
$ \ok $ be an algebraic closure of
$ k $, and put
$ Y_{\ok}^{(i)} = Y^{(i)}\otimes_{k}\ok $.
Let
$ l $ be a prime number different from the characteristic of
$ k $.
Then the
$ E_{1} $-term of the weight spectral sequence of Rapoport and Zink
is given by direct sums of the
$ l $-adic cohomology groups of
$ Y_{\ok}^{(i)} $ which are Tate-twisted appropriately.
Its differential
$ d_{1} $ consists of the Cech restriction morphisms and the
co-Cech Gysin morphisms which are denoted respectively by
$ \rho $ and
$ \gamma $ in this paper.
The primitive cohomology of
$ Y_{\ok}^{(i)} $ has a canonical pairing induced by
Poincar\'e duality together with the hard Lefschetz theorem [4]
(choosing and fixing an ample divisor
class of
$ f) $.

\medskip\noindent
{\bf 0.2.~Theorem.}
{\it Conjecture {\rm (0.1)} is true if the restrictions of the
canonical pairing to the intersections of the primitive part with
$ \Im\,\rho $ and with
$ \Im\,\gamma $ are both nondegenerate.
}

\medskip
So the problem is reduced to the study of the canonical pairing on
the primitive part.
There are some examples satisfying the assumption of (0.2),
see (2.8).
If
$ n = 3 $, the converse of (0.2) is also true and the hypothesis on
$ \Im\,\gamma $ is always satisfied in this case (using [2], [21]).
Note that (0.1) may apparently depend on the choice of
$ l $, and we can prove (0.1) for certain
$ l $ in a simple case where every eigenvalue of the Frobenius
action has multiplicity
$ 1 $, see (2.7) and (5.4).
(In this case we can determine the endomorphism ring of a simple
motive, see (5.2).)
However, the general case seems to be related closely to the
standard conjectures.

A conjecture similar to (0.2) but using only the restriction
morphisms was noted by Rapoport and Zink in the introduction of
[21].
We can prove this conjecture under an inductive hypothesis as
follows:

\medskip\noindent
{\bf 0.3.~Theorem.}
{\it Assume that Conjecture {\rm (0.1)} holds for a general
hyperplane section of the generic fiber, and that the
restriction of the modified pairing in {\rm [21]} to
$ \Im\,\rho $ or the restriction of the canonical pairing to
the intersection of the primitive part with
$ \Im\,\rho $ is nondegenerate.
Then {\rm (0.1)} is true.
}

\medskip
In the complex analytic case, the hypothesis of (0.2) is trivially
satisfied because of the positivity of polarizations of Hodge structures.
We can argue similarly if we have a kind of ``positivity'' in
characteristic
$ p > 0 $.
The positivity for a zero-dimensional variety (or a motive) is clear,
because the pairing is defined over the subfield
$ \bQ $ of
$ \bQ_{l}) $.
In the one-dimensional case, this notion is provided by the theory
of Riemann forms for abelian varieties [31] (see also [17], [20])
combined with work of Deligne [6] on the compatibility of the
Weil pairing and Poincar\'e duality (see also (3.4) below).
These were used in an essential way for the proof of Conjecture (0.1)
for
$ n = 2 $ in [21].
However, for higher dimensional varieties, we do not know any notion
of positivity except the standard conjecture of Hodge index type [14].
Using the theory on the standard conjectures (loc. cit.)
we show

\medskip\noindent
{\bf 0.4.~Theorem.}
{\it Assume that the standard conjectures
hold for
$ Y^{(i)} $,
$ Y^{(i)}\mtim_{k}Y^{(i)} $, and the numerical equivalence and the
homological equivalence coincide for
$ Y^{(i)}\mtim_{k}Y^{(i+1)} \,(i > 0) $.
Then Conjecture {\rm (0.1)} is true.
}

\medskip
I am quite recently informed that if the generic fiber can be
uniformized by the Drinfeld upper half space, Conjecture (0.1)
is proved by T. Ito using Theorem (0.2).

Part of this work was done during my stay at the University of Leiden.
I thank Professor J. Murre for useful discussions, and
the staff of the institute for the hospitality.
I also thank Takeshi Saito for a useful discussion about the
consequence of de Jong's theory of alternations.

In Sect.~1 and Sect.~2, we review the theory of graded or bigraded
modules of Lefschetz-type and prove (0.2) and (0.3).
In Sect.~3, we review the work of Deligne on the compatibility of
the Weil pairing and Poincar\'e duality.
In Sect.~4, we review the standard conjectures and prove (0.4).
In Sect.~5, we study the Frobenius action, and prove (0.1) in some
simple cases.

\bigskip\bigskip
\centerline{{\bf 1. Graded Modules of Lefschetz-Type}}

\bigskip\noindent
We first review a theory of morphisms of degree
$ 1 $ between graded modules of Lefschetz-type [23].
A typical example is given by the restriction or Gysin morphism
associated to a morphism of smooth projective varieties whose
relative dimension is
$ -1 $.

\medskip\noindent
{\bf 1.1.}
In this and the next sections we denote by
$ \Lambda $ a field.
By a graded
$ \Lambda[L] $-module we will mean a finite dimensional
graded vector space
$ M^{\ssbull} $ over
$ \Lambda $ having an action of
$ L $ with degree
$ 2 $.
We call
$ M^{\ssbull} $
$ n $-{\it symmetric} if
$$
L^{j} : M^{n-j} \simto M^{n+j} \quad\text{for}\,\, j > 0.
$$
Here the Tate twists are omitted to simplify the notation.
We say that
$ M^{\ssbull} $ is a graded module of Lefschetz-type if it is a
$ 0 $-symmetric graded
$ \Lambda[L] $-module.
Then we have the Lefschetz decomposition
$$
M^{j} = \mopls_{i\ge 0}\, L^{i}{}_{0}M^{j-2i}\quad\text{with}
\,\,{}_{0}M^{-j} = \Ker\,L^{j+1} \subset M^{-j}.
$$

We say that
$ f : M^{\ssbull} \to N^{\ssbull} $ is a morphism of degree
$ m $ between graded modules of Lefschetz-type, if
$ f(M^{j}) \subset N^{j+m} $ and
$ f $ is
$ \Lambda[L] $-linear.
Shifting the degree of
$ N^{\ssbull} $, it is identified with a graded morphism of a
$ 0 $-symmetric graded
$ \Lambda[L] $-module to a
$ (-m) $-symmetric graded
$ \Lambda[L] $-module.
Considering the image of the primitive part, we see that a
morphism of degree
$ 0 $ preserves the Lefschetz decomposition and there is no
nontrivial morphism of negative degree (i.e. a morphism of a
$ 0 $-symmetric module to an
$ n $-symmetric module for
$ n > 0 $ is trivial).
For a morphism of degree
$ 1 $, we see that
$$
f({}_{0}M^{-j}) \subset {}_{0}N^{-j+1} + L\,{}_{0}N^{-j-1}.
\leqno(1.1.1)
$$

Let
$ f : M^{\ssbull} \to N^{\ssbull} $ be a morphism of degree
$ 1 $ between graded modules of Lefschetz-type.
We define
$$
\Im^{0}f = \mopls_{j\in\bZ}(\mopls_{i\ge 0}\,
L^{i}(\Im\,f \cap {}_{0}N^{-j})),
\quad \Im^{1}f = \Im\,f/\Im^{0}f.
$$
Then
$ \Im^{0}f $ is
$ 0 $-symmetric, and its primitive part is given by
$ \Im\,f \cap {}_{0}N^{-j} $.
It has been remarked by T. Ito that the proof of the next lemma
easily follows from the definition itself.

\medskip\noindent
{\bf 1.2.~Lemma.}
{\it The quotient module
$ \Im^{1}f $ is
$ 1 $-symmetric, i.e.
we have the bijectivity of
$$
L^{j+1} : (\Im^{1}f)^{-j} \to (\Im^{1}f)^{j+2}\quad
\text{for}\,\, j \ge 0,
\leqno(1.2.1)
$$
where the upper index
$ -j $ means the degree
$ -j $ part, etc.
}

\medskip\noindent
{\it Proof.}
The surjectivity follows from the
$ 0 $-symmetry of
$ M $.
To show the injectivity, let
$ n = \sum_{i\ge 0}\,L^{i}n_{i} \in (\Im\,f)^{-j} $ with
$ n_{i} \in {}_{0}N^{-j-2i} $, and assume
$ L^{j+1}n \in (\Im^{0}f)^{j+2} $.
Since
$ L^{j+1}n = \sum_{i>0}\,L^{j+i+1}n_{i} $ and
$ \Im^{0}f $ is
$ 0 $-symmetric, we may assume
$ n_{i} = 0 $ for
$ i > 0 $ by modifying
$ n $ modulo
$ \Im^{0}f $ if necessary.
Then
$ n \in \Im^{0}f $ by definition, and the assertion follows.

\medskip\noindent
{\bf 1.3.~Proposition.}
{\it Let
$ M^{\ssbull}, N^{\ssbull} $ be graded
$ \Lambda[L] $-modules of Lefschetz-type, and
$$
f : M^{\ssbull} \to N^{\ssbull},\quad
g : N^{\ssbull} \to M^{\ssbull}
$$
be morphisms of degree one.
Assume there are nondegenerate pairings of
$ \Lambda $-modules
$$
\Phi_{M} : M^{j}\otimes_{\Lambda}M^{-j} \to \Lambda,\quad
\Phi_{N} : N^{j}\otimes_{\Lambda}N^{-j} \to \Lambda,
$$
such that
$ \Phi_{N}\scirc (f\otimes id) $ coincides with
$ \Phi_{M}\scirc (id\otimes g) $ up to a nonzero multiple constant
and that
$ \Phi_{M}\scirc (id\otimes L) = \Phi_{M}\scirc (L\otimes id) $
{\rm (}and the same for
$ \Phi_{M}) $.
Then
}
$$
\dim\,(\Im^{0}f)^{j} = \dim\,(\Im^{1}g)^{j+1},\quad
\dim\,(\Im^{0}g)^{j} = \dim\,(\Im^{1}f)^{j+1}.
$$

\medskip\noindent
{\it Proof.}
Let
$ a_{j}, b_{j}, c_{j}, d_{j} $ denote respectively the above
dimensions so that
$ a_{j} = a_{-j} $,
$ b_{j} = b_{-j} $, etc.
By the duality of
$ f $ and
$ g $, we have
$$
a_{j+1} + d_{j} = \dim\,(\Im\,f)^{j+1} = \dim\,(\Im\,g)^{-j} =
b_{-j-1} + c_{-j} = b_{j+1} + c_{j}.
$$
If we put
$ p_{j} = a_{j} - b_{j} $,
$ q_{j} = c_{j} - d_{j} $, we get
$ p_{j+1} = q_{j} $ for any
$ j \in \bZ $.
So
$ p_{j} = q_{j} = 0 $ by the symmetry of
$ p_{j} $,
$ q_{j} $, and the assertion follows.

\medskip\noindent
{\bf 1.4.~Lemma.}
{\it With the notation and assumption of {\rm (1.3),} assume that
$ n \in \Im\, f \cap {}_{0}N^{-j} $ vanishes if
$ \Phi_{N}(f(m),L^{j}n) = 0 $ for any
$ m \in M^{-j-1} $.
Then
}
$$
\Ker\,g \cap \Im^{0}f = 0.
\leqno(1.4.1)
$$

\medskip\noindent
{\it Proof.}
Let
$ n = \sum_{0\le i\le r} L^{i}n_{i} $ with
$ n_{i} \in {}_{0}N^{-j-2i} $,
and assume
$ g(n) = 0 $.
Since
$$
L^{j+2r}g(n_{r}) = L^{j+r}g(n) = 0,
$$
the assertion is reduced by induction on
$ r $ to the injectivity of
$$
g\scirc L^{j} : \Im\, f \cap{}_{0}N^{-j} \to M^{j+1}.
$$
But this injectivity is clear, because
$ \Phi_{N}(f(m),L^{j}n) $ is equal to
$ \Phi_{M}(m,g(L^{j}n)) $ up to a nonzero multiple constant.
So the assertion follows.

\medskip\noindent
{\bf 1.5.~Proposition.}
{\it With the notation and assumption of {\rm (1.3),} assume
the vanishing of
$ fgf $
and the injectivity of
$$
f : \Im^{0}g \to N^{\ssbull},\quad
g : \Im^{0}f \to M^{\ssbull}.
\leqno(1.5.1)
$$
Then the compositions
$$
\Im^{0}g \overset f\to\to \Im\,f \to
\Im^{1}f,\quad \Im^{0}f
\overset g\to\to \Im\,g \to \Im^{1}g
\leqno(1.5.2)
$$
are isomorphisms, and we have canonical decompositions
$$
\Im\,f = \Im^{0}f\oplus \Im^{1}f,\quad
\Im\,g = \Im^{0}g\oplus \Im^{1}g,
\leqno(1.5.3)
$$
such that the restriction of
$ g $ to
$ \Im^{0}f $ is injective, and that to
$ \Im^{1}f $ is zero, and similarly for the restriction of
$ f $.
Furthermore, we have canonical isomorphisms
}
$$
\Im\,fg = \Im^{1}f,\quad
\Im\,gf = \Im^{1}g.
\leqno(1.5.4)
$$

\medskip\noindent
{\it Proof.}
The hypothesis implies
$ f(\Im^{0}g) \cap \Im^{0}f = g(\Im^{0}f) \cap
\Im^{0}g = 0 $, because the vanishing of
$ fgf $ is equivalent to that of
$ gfg $ by duality.
So the assertion follows from (1.3).

\medskip\noindent
{\bf 1.6.~Proposition.}
{\it With the notation and the assumptions of {\rm (1.3),} assume
further that {\rm (1.4.1)} holds,
$ \Im\, fg $ is
$ 1 $-symmetric,
$ gfg = 0 $, and
$$
(\Ker\,g \cap \Im\,f)^{j} = (\Im\,fg)^{j}\quad \text{for}\,\,
j \ne 0.
\leqno(1.6.1)
$$
Then {\rm (1.6.1)} holds also for
$ j = 0 $.
}

\medskip\noindent
{\it Proof.}
Since
$ \Ker\,g \cap \Im\,f \supset \Im\,fg $,
it is enough to show the coincidence of the images of both sides
of (1.6.1) in
$ \Im^{1}f $ using (1.4.1).
We can identify
$$
(\Ker\,g \cap \Im\,f)/\Im\,fg
\leqno(1.6.2)
$$
with a graded
$ \Lambda [L] $-submodule of
$ \Im^{1}f/\Im\,fg $, and the last module is
$ 1 $-symmetric by hypothesis and (1.2).
Furthermore (1.6.2) vanishes except for the degree
$ 0 $ by hypothesis.
So it vanishes at any degree, and the assertion follows.

\newpage
\centerline{{\bf 2. Bigraded Modules of Lefschetz-Type}}

\bigskip\noindent
We prove (0.2--3) after reviewing a theory of bigraded modules of
Lefschetz-type [23], [24].
These modules appear in the
$ E_{1} $-term of the Steenbrink-type spectral sequence
associated to a semi-stable degeneration.

\medskip\noindent
{\bf 2.1.}
Let
$ M^{\ssbull ,\ssbull} $ be a bigraded
$ \Lambda[N,L] $-module of Lefschetz-type, i.e.
it is a finite dimensional bigraded vector space over
$ \Lambda $ having commuting actions of
$ N, L $ with bidegrees
$ (2,0) $ and
$ (0.2) $ respectively such that
$$
N^{i} : M^{-i,j} \simto M^{i,j} \,(i > 0),\quad L^{j} :
M^{i,-j} \simto M^{i,j} \,(j > 0).
$$
Put
$ M_{i}^{j} = M^{-i,j} $, and
$ {}_{(0)}M_{i}^{j} = \Ker\,N^{i+1} \subset M_{i}^{j} $ for
$ i \ge 0 $,
and
$ 0 $ otherwise.
Then we have the Lefschetz decomposition for the first index:
$$
M_{i}^{j} = \mopls_{a\ge 0}\, N^{a} {}_{(0)}M_{i+2a}^{j}.
\leqno(2.1.1)
$$
We define
$ C_{i}^{j} = {}_{(0)}M_{i}^{j} $ and
$$
C_{i,a}^{j} = N^{a} C_{i+a}^{j} \subset M_{i-a}^{j},
$$
so that
$ C_{i,a}^{j} = 0 $ for
$ i < 0 $ or
$ a < 0 $,
and
$ C_{i,0}^{j} = C_{i}^{j} $.
Then
$$
N^{i-a} : C_{i,a}^{j} \simto C_{a,i}^{j} \,(i > a),\quad
L^{j} : C_{i,a}^{-j} \simto C_{i,a}^{j} \,(j > 0).
\leqno(2.1.2)
$$

Let
$ d $ be a differential of bidegree
$ (1,1) $ on
$ M^{\ssbull ,\ssbull} $ which commutes with
$ N, L $ and satisfies
$ d^{2} = 0 $.
By (1.1.1) applied to the action of
$ N $, we have a decomposition
$ d = d' + d'' $ such that
$$
d' : C_{i,a}^{j} \to C_{i-1,a}^{j+1},\quad
d'' : C_{i,a}^{j} \to C_{i,a+1}^{j+1}
$$
are differentials which anti-commute with each other.
Let
$ C_{i}^{\ssbull} = \mopls_{j\in\bZ}C_{i}^{j} $.
We have morphisms of degree
$ 1 $ between graded
$ \Lambda[L] $-modules
$$
\gamma_{i} : C_{i}^{\ssbull} \to C_{i-1}^{\ssbull}\,\,(i > 0),
\quad
\rho_{i} : C_{i}^{\ssbull} \to C_{i+1}^{\ssbull}\,\,(i \ge 0)
$$
such that
$ d', d'' $ are identified with
$ \gamma_{i-a}, \rho_{i-a} $ respectively.

We define
$ H_{i}^{j} = Z_{i}^{j}/B_{i}^{j} $ with
$$
Z_{i}^{j} = \Ker(d : M_{i}^{j} \to M_{i-1}^{j+1}),\quad
B_{i}^{j} = \Im(d : M_{i+1}^{j-1} \to M_{i}^{j}).
$$
Then we have the induced morphism
$ N : H_{i}^{j} \to H_{i-2}^{j} $, and similarly for
$ Z_{i}^{j}, B_{i}^{j} $.
For a positive integer
$ i $ and an integer
$ j $, we will consider the condition for the bijectivity of
the morphism
$$
N^{i} : H_{i}^{j} \to H_{-i}^{j}\quad\text{for}\,\,i > 0.
\leqno(2.1.3)
$$
For this we will assume that the action of
$ L $ induces bijections
$$
L^{j} : H_{i}^{-j} \simto H_{i}^{j}\quad\text{for}\,\,j > 0,
\leqno(2.1.4)
$$
and that both terms of (2.1.3) have the same dimension (so that
bijectivity is equivalent to injectivity and to surjectivity).
Using (2.1.4), the last assumption is satisfied if we have a
self-duality of
$ M^{\ssbull,\ssbull} $.

Since
$ M_{i-1}^{\ssbull} = C_{i-1}^{\ssbull} \oplus
NM_{i+1}^{\ssbull} $ and
$ N^{i} : M_{i}^{j} \simto M_{-i}^{j} \,(i > 0) $,
we have a morphism
$$
\td:= d\scirc N^{-i} : Z_{-i}^{j} \to
C_{i-1}^{j+1}\quad\text{for}\,\,i > 0.
$$

\medskip\noindent
{\bf 2.2.~Proposition.}
{\it {\rm (i)} The surjectivity of {\rm (2.1.3)} is equivalent
to
$$
\td(Z_{-i}^{j}) = (\Im\,\gamma_{i}\rho_{i-1})^{j+1}.
\leqno(2.2.1)
$$
{\rm (ii)} If {\rm (2.1.3)} for
$ (i + 2, j) $ is surjective, then the surjectivity of
{\rm (2.1.3)} is further equivalent to
$$
\gamma_{i}(\Ker\,\rho_{i})^{j} =
(\Im\,\gamma_{i}\rho_{i-1})^{j+1}.
\leqno(2.2.2)
$$
{\rm (iii)} If furthermore {\rm (2.1.3)} for
$ (i + 2, j) $ is surjective and {\rm (2.1.3)} for
$ (i + 1, j + 1) $ is injective, then the surjectivity of
{\rm (2.1.3)} is further equivalent to
$$
(\Ker\,\rho_{i-1} \cap \Im\,\gamma_{i})^{j+1} =
(\Im\,\gamma_{i}\rho_{i-1})^{j+1}.
\leqno(2.2.3)
$$
{\rm (iv)} If {\rm (2.1.3)} for
$ (i + 1, j - 1) $ is surjective, then the injectivity of
{\rm (2.1.3)} is equivalent to
}
$$
(\Ker\,\gamma_{i} \cap \Im\,\rho_{i-1})^{j} =
(\Im\,\rho_{i-1}\gamma_{i})^{j}\,
(= (\Im\,\gamma_{i+1}\rho_{i})^{j}).
\leqno(2.2.4)
$$

\medskip\noindent
{\it Proof.}
We have by definition
$ \td(B_{-i}^{j}) = (\Im\,\gamma_{i}\rho_{i-1})^{j+1} $.
Then the first assertion is clear because the surjectivity is
equivalent to
$ Z_{-i}^{j} = B_{-i}^{j} + N^{i}Z_{i}^{j} $.

For the second, take
$ m = \sum_{a\ge 0}\,N^{a}m_{a} \in N^{-i}Z_{-i}^{j} $ with
$ m_{a} \in C_{i+2a}^{j} $ to calculate
$ \td(Z_{-i}^{j}) $.
We may assume that
$ m_{0} \in \Ker\,\rho_{i} $ applying (2.2.1) for
$ i + 2 $ to
$ \sum_{a\ge 1}N^{a-1}m_{a} $ and modifying
$ m $ by an element of
$ d(C_{i+1}^{j-1}) $ because this does not change
$ dm $.
So the second assertion is clear.

Before showing the third assertion, we see that
$$
(\Ker\,\gamma_{i} \cap \Im\,\rho_{i-1})^{j}/\td(Z_{-i-1}^{j-1})
\simto (N^{-i}B_{-i}^{j} \cap Z_{i}^{j})/B_{i}^{j},
$$
because
$ N^{-i}B_{-i}^{j} = (\Im\,\rho_{i-1})^{j} + B_{i}^{j} $.
Then the last assertion follows from (2.2.1).

For the third assertion, take
$ \gamma_{i}m \in (\Ker\,\rho_{i-1})^{j+1} $, where
$ m \in C_{i}^{j} $ with
$ \rho_{i-1}\gamma_{i}m = 0 $.
Then
$ \rho_{i}m \in \Ker\,\gamma_{i+1} $,
and we may assume
$ \rho_{i}m = 0 $ using (2.2.4) for
$ (i+1,j+1) $ and modifying
$ m $ by
$ \gamma_{i-1}m' $ because it does not change
$ \gamma_{i}m $.
So the assertion follows from (2.2.2).

\medskip\noindent
{\bf 2.3.~Application.}
With the notation of the introduction we may assume
$$
C_{i}^{j} = H^{n-i+j}(Y_{\ok}^{(i+1)},\bQ_{l}(-i)),\quad
H_{i}^{j} = \Gr_{n+j+i}^{W}H^{n+j}(X_{\oK},\bQ_{l}),
\leqno(2.3.1)
$$
and
$ \rho_{i}, \gamma_{i+1} $ are respectively the Cech restriction
and co-Cech Gysin morphisms, which are dual of each other up to
a sign.
In particular,
$ \rho_{j}\gamma_{j+1}\rho_{j} = -
\gamma_{j+2}\rho_{j+1}\rho_{j} = 0 $, i.e. the first hypothesis
of Proposition (1.5) is satisfied.
Note that
$ N $ is the logarithm of the monodromy, and
$ L $ is given by the ample divisor class of
$ f $, see [11], [21].
So the bijectivity of (2.1.4) follows from the hard Lefschetz
theorem for the generic fiber together with the strict
compatibility of the weight filtration [4].

\medskip\noindent
{\bf 2.4.~Conjecture of Rapoport and Zink.}
For a smooth projective
$ \ok $-variety
$ Y $ with an ample line bundle
$ L $, we have a canonical pairing on
$ H^{j}(Y,\bQ_{l}) $ by Poincar\'e duality and the Hard
Lefschetz theorem [4].
Using further the Lefschetz decomposition and modifying the sign
as in (4.1), we get a modified pairing on
$ H^{j}(Y,\bQ_{l}) $ as in [21].
Rapoport and Zink noted there that Conjecture (0.1) would be
verified if the restriction of this modified pairing to
$ \Im\,\rho $ is nondegenerate.
Note that the hypothesis of (1.4) is satisfied under the above
hypothesis (where
$ f = \rho_{i-1} $ and
$ g = \gamma_{i} $) because the modified pairing coincides with
the canonical pairing up to a sign if one factor belongs to the
primitive part.

\medskip\noindent
{\bf 2.5.~Proof of (0.2) and (0.3).}
By [21], the
$ E_{1} $-term of the weight spectral sequence has a
structure of bigraded
$ \bQ_{l} $-modules of Lefschetz-type, see (2.3).
Then (2.2.3--4) follows from (1.4--5), and Theorem (0.2) follows.

For Theorem (0.3) we will show by decreasing induction on
$ i $ that (2.1.3) is bijective (or equivalently, injective) and
$ \Im\,\gamma_{i}\rho_{i-1} $ is
$ 1 $-symmetric.
Assume the two assertions are true for
$ i + r $ with
$ r > 0 $.
Considering the restriction morphism to a general hyperplane
section of the generic fiber and using the weak Lefschetz
theorem, we see that (2.1.3) is injective for
$ j < 0 $,
because the restriction morphism is strictly compatible with the
weight filtration.
Then (2.1.3) is injective also for
$ j > 0 $ by (2.1.4), and it is injective for any
$ j $ by (2.2) (iv) and (1.6) where
$ f = \rho_{i-1} $ and
$ g = \gamma_{i} $.
So it remains to show that
$ \Im\,\gamma_{i}\rho_{i-1} $ is
$ 1 $-symmetric.

Since the injectivity of the first morphism of (1.5.1)
follows from (1.4), the assertion is reduced by (1.5) to the
injectivity of
$ \rho_{i-1} : \Im^{0}\gamma_{i} \to C_{i}^{\ssbull}[1] $,
i.e. to the vanishing of
$$
\Ker(\rho_{i-1} : \Im^{0}\gamma_{i} \to
(\Im\,\rho_{i-1}\gamma_{i})[1]).
\leqno(2.5.1)
$$
Here
$ [m] $ means a shift of degree for an integer
$ m $, i.e.
$ (M[m])^{i} = M^{i+m} $.
We see that (2.5.1) is
$ 0 $-symmetric, because
$ \Im^{0}\gamma_{i} $ and
$ (\Im\,\rho_{i-1}\gamma_{i})[1] $ are (using inductive hypothesis).
Furthermore (2.5.1) is a submodule of
$ \Ker\,\rho_{i-1} \cap \Im\,\gamma_{i} $, and the latter is
identified by (2.2.3--4) with
$$
\Im(\gamma_{i} : (\Im\,\rho_{i-1})[-1] \to \Im\,\gamma_{i}) =
(\Im\,\rho_{i-1}/\Im\,\rho_{i-1}\gamma_{i})[-1].
$$
The last module is an extension of a
$ 2 $-symmetric module by a
$ 1 $-symmetric module, because
$ \Im\,\rho_{i-1}\gamma_{i}\cap\Im^{0}\rho_{i-1} = 0 $ by (1.4)
and
$ \Im\,\rho_{i-1}\gamma_{i}\,(= \Im\,\gamma_{i+1}\rho_{i}) $ is
$ 1 $-symmetric by inductive hypothesis.
But there is no nontrivial morphism of a
$ 0 $-symmetric module to an
$ n $-symmetric module for
$ n > 0 $, see (1.1).
So the assertion follows.

\medskip\noindent
{\bf 2.6.~First cohomology case.}
The restriction of the pairing to the image of
$ \rho_{i-1} $ in
$ H^{1}(Y_{\ok}^{(i+1)},\bQ_{l}) $ is always nondegenerate by
the abelian-positivity in Sect.~3, because
$ \rho_{i-1} $ is associated to the morphism of Picard varieties
using the Tate module.

\medskip\noindent
{\bf 2.7.~Three-dimensional case.} Assume
$ n = 3 $.
Then we can show that Conjecture (0.1) is equivalent to the
nondegeneracy of the restrictions of the canonical pairing to
$$
H^{3}(Y_{\ok}^{(1)},\bQ_{l})^{\prim} \cap \Im\,\gamma,\quad
H^{2}(Y_{\ok}^{(2)},\bQ_{l})^{\prim} \cap \Im\,\rho.
\leqno(2.7.1)
$$
Indeed, using the abelian-positivity in Sect.~3, the arguments
in (2.5) show that these conditions in this case are equivalent
respectively to the bijectivity of
$$
N : H_{1}^{-1} \to H_{-1}^{-1},\quad
N : H_{1}^{0} \to H_{-1}^{0}.
$$
Since the first morphism is bijective by using a general
hyperplane section, the nondegeneracy of the pairing for
the intersection with
$ \Im\,\gamma $ in (0.2) is always true for
$ n = 3 $.
Note that the bijectivity of the first morphism cannot be proved by
using simply the abelian positivity in Sect.~3 unless we know that
$ \Ker\,\gamma \cap H^{1}(Y_{\ok}^{(2)},\bQ_{l}) $ corresponds to
an abelian subvariety of the Picard variety.

If
$ \Im\,\rho \subset H^{2}(Y_{\ok}^{(2)},\bQ_{l}) $ is contained
in the subspace generated by algebraic cycle classes, then
Conjecture (0.1) can be reduced to
$ D(Y^{(1)}) $ or
$ A(Y^{(1)},L) $ in the notation of (4.1)
(or to the Tate conjecture [30] for divisors on
$ Y^{(1)} $) because
$ D(Y^{(2)}) $,
$ D(Y^{(1)}) $ for divisors and
$ I(Y^{(2)},L) $ are known.

If, for every eigenvalue of the Frobenius action on
$ \Im\, \rho \cap H^{2}({Y}_{\ok}^{(2)},\bQ_{l})^{\prim} $, its
multiplicity as an eigenvalue of the Frobenius action on
$ H^{2}({Y}_{\ok}^{(2)},\bQ_{l})^{\prim} $ is
$ 1 $, then Conjecture (0.1) can be proved for certain prime
numbers
$ l $, see (5.4).

\medskip\noindent
{\bf 2.8.~Example.}
Let
$ V \to S := \Spec R $ be a smooth projective morphism of
relative dimension
$ n + 1 $.
For
$ i = 0, \dots, r $,
let
$ Z_{i} $ be a smooth hypersurface of
$ V $ defined by a section
$ P_{i} $ of a relative ample line bundle
$ L_{i} $.
Assume
$ L_{0} = \bigotimes_{1\le i\le r} L_{i} $, and
$ \bigcup_{0\le i\le r} Z_{i} $ is a divisor with normal crossings
relative to
$ R $ (i.e. the fiber over
$ \Spec k $ is a divisor with normal crossings).
Let
$ \pi $ be a generator of the maximal ideal of
$ R $, and define
$$
P = P_{1}\cdots P_{r} + P_{0}\pi.
$$
Let
$ X $ be the hypersurface of
$ V $ defined by
$ P $.
(This is an analogue of [7] for
$ r = 2 $.)
Then Conjecture (0.1) is true for a semi-stable model of the
generic fiber of
$ X $ although
$ X $ is not semistable in general.

Indeed, consider an iteration of blow-ups
$ \sigma_{j} : V^{(j)} \to V^{(j-1)} $ along the proper
transform of
$ Z_{0} \cap Z_{j} $ for
$ j = 1, \dots, r - 1 $,
where
$ V^{(0)} = V $.
We have local coordinates
$ x_{0}, \dots, x_{n} $ over
$ R $ such that
$ P $ is locally given by
$$
ux_{1} \cdots x_{m} + x_{0}\pi,
$$
where
$ x_{i} $ is the restriction of
$ P_{i} $ for
$ i \le m $ and
$ u $ is locally invertible.
Then the proper transform of
$ P $ by the blow-up along
$ Z_{0} \cap Z_{1} $ is locally given by
$$
ux_{2} \cdots x_{m} + x'_{0}\pi \quad \text{or}\quad
ux''_{1}x_{2} \cdots x_{m} + \pi,
$$
where
$ (x'_{0}, x_{1}, x_{2}, \cdots, x_{n}) $ and
$ (x_{0}, x''_{1}, x_{2}, \cdots, x_{n}) $ are local
coordinate systems on open subvarieties
$ U' $,
$ U'' $ of
$ V^{(1)} $ such that
$ x_{0} = x'_{0}x_{1} $ on
$ U' $ and
$ x_{1} = x''_{1}x_{0} $ on
$ U'' $.
Here the pull-back of
$ x_{i} $ is also denoted by
$ x_{i} $ to simplify the notation.
Since
$ U'' $ does not intersect the proper transform of
$ Z_{0}\cap Z_{j} $ for
$ j > 1 $,
we can proceed by induction, and get a semi-stable model.
Its generic fiber is same as that of
$ X $ because the intersection of the center of the blow-up with
the generic fiber of (the proper transform of)
$ X $ is a locally principal divisor on the generic fiber.
Furthermore the
$ Y^{(i)} $ are lifted to smooth projective schemes over
$ R $ so that the assumption of (0.2) is verified by using
Hodge theory.

\bigskip\bigskip
\centerline{{\bf 3. Abelian-Positivity}}

\bigskip\noindent
{\bf 3.1.~Canonical pairing.}
Let
$ A $ be an abelian variety over a field
$ k $.
We denote by
$ A^{\vee} $ its dual variety, and by
$ T_{l}A_{\ok} $ the Tate module of
$ A_{\ok} := A\otimes_{k}\ok $ where
$ l $ is a prime number different form
$ \char k $.
Using the Kummer sequence, we have a canonical isomorphism
$$
H^{1}(A_{\ok},\mu_{n}) = A^{\vee}(\ok)_{n} \,
(:= \Ker(n : A^{\vee}(\ok) \to A^{\vee}(\ok))),
\leqno(3.1.1)
$$
where
$ n $ is an integer prime to
$ \char k $.
Then, passing to the limit, we get
$$
H^{1}(A_{\ok},\bZ_{l}(1)) = T_{l}{A}_{\ok}^{\vee}.
\leqno(3.1.2)
$$
Since the left-hand side of (3.1.1) is identified with
$$
\Hom(A(\ok)_{n}, \mu_{n}),
$$
using torsors [6], we get the canonical pairing of Weil [31]
(see also [17], [20]):
$$
A(\ok)_{n}\otimes A^{\vee}(\ok)_{n} \to
\mu_{n},\quad T_{l}A_{\ok}
\otimes_{\bZ_{l}}T_{l}{A}_{\ok}^{\vee} \to \bZ_{l}(1).
\leqno(3.1.3)
$$
To get a pairing of
$ T_{l}A_{\ok} $,
we take a divisor
$ D $ on
$ A $ which induces a morphism
$$
\varphi_{D} : A \to A^{\vee},
\leqno(3.1.4)
$$
such that
$ \varphi_{D}(a) \in A^{\vee}(\ok) $ for
$ a \in A(\ok) $ is given by
$ {T}_{a}^{*}D_{\ok} - D_{\ok} $,
where
$ T_{a} $ is the translation by
$ a $ (see loc. cit.)
Note that
$ \varphi_{D} $ depends only on the algebraic equivalence class of
$ D $,
and
$ \varphi_{D} $ is an isogeny if
$ D $ is ample.

If the pairing is induced by an ample divisor on
$ A $,
its restriction to
$ T_{l}B_{\ok} $ for any abelian subvariety
$ B $ of
$ A $ is nondegenerate, because we have the commutative diagram:
$$
\CD
B @>>> A
\\
@VV{\varphi_{D|_{B}}}V @VV{\varphi_{D}}V
\\
B^{\vee} @<<< A^{\vee}
\endCD
\leqno(3.1.5)
$$
Note that this holds only for subgroups of
$ T_{l}A_{\ok} $ corresponding to abelian subvarieties.

We say that a pairing of a
$ \bQ_{l} $-module
$ V $ with a continuous action of
$ G := \Gal(\ok/k) $ is abelian-positive if there exists an
abelian variety with an ample divisor
$ D $ such that
$ V $ is isomorphic to
$ T_{l}A_{\ok}\otimes_{\bZ_{l}}\bQ_{l} $ up to a Tate
twist as a
$ \bQ_{l}[G] $-module and the pairing corresponds to the one on
$ T_{l}A_{\ok} $ defined by the canonical pairing and
$ \varphi_{D} $.
Note that abelian-positive pairings (having the same weight) are
stable by direct sums.

\medskip\noindent
{\bf 3.2.~Compatibility of the cycle classes.}
Let
$ A $ be an abelian variety over
$ k $,
$ X $ a smooth projective variety over
$ k $,
and
$ D $ a divisor on
$ A \mtim_{k} X $ such that its restriction to
$ \{0\} \mtim X $ is rationally equivalent to zero.
Let
$ P $ be the Picard variety of
$ X $.
Then
$ D $ induces a morphism of abelian varieties
$$
\Psi_{D} : A \to P,
$$
such
$ \Psi_{D}(a) \in P(\ok) $ for
$ a \in A(\ok) $ is defined by the restriction of
$ D_{\ok} $ to
$ \{a\} \mtim X_{\ok} $,
see [31] (and also [17], [20]).

Let
$ cl(D)^{1,1} \in H^{1}(A_{\ok}, R^{1}(pr_{1})_{*}\mu_{n}) $
denote the
$ (1,1) $-component of the cycle class of
$ D $,
where
$ n $ is an integer prime to
$ \char k $, and
$ pr_{1} $ is the first projection.
Assume
$ \NS(X) $ is torsion-free.
Then
$ R^{1}(pr_{1})_{*}\mu_{n} $ is a constant sheaf on
$ A_{\ok} $ with fiber
$ H^{1}(X_{\ok},\mu_{n}) = P(\ok)_{n} $,
and we get
$$
cl(D)^{1,1} \in H^{1}(A_{\ok}, R^{1}(pr_{1})_{*}\mu_{n}) = \Hom
(A(\ok)_{n},P(\ok)_{n}).
$$

\medskip\noindent
{\bf 3.3.~Theorem} {\rm (Deligne).}
{\it The induced morphism
$ \Psi_{D} : A(\ok)_{n} \to P(\ok)_{n} $ coincides with
$ - cl(D)^{1,1} $.
}

\medskip
(The proof is essentially the same as in [6].)

\medskip\noindent
{\bf 3.4.~Corollary} {\rm (Deligne).}
{\it Let
$ C $ be a smooth projective curve over a field
$ k $, and
$ J $ its Jacobian.
Then we have a canonical isomorphism
$ H^{1}(C_{\ok},\bZ_{l}(1)) = T_{l}J(\ok) $ such that
Poincar\'e duality on
$ H^{1}(C_{\ok},\bZ_{l}(1)) $ is identified with the pairing of
$ T_{l}J(\ok) $ given by the canonical pairing {\rm (3.1.3)}
together with {\rm (3.1.4)} for the theta divisor on
$ J $.
}

\medskip\noindent
{\it Proof.}
We may assume that
$ C $ has a
$ k $-rational point replacing
$ k $ with a finite extension
$ k' $ and
$ C $ with
$ C \otimes_{k}k' $ if necessary.
Choosing a
$ k $-rational point of
$ C $,
we have a morphism
$ f : C \to J $.
It is well-known that this induces isomorphisms
$$
\aligned
&f^{*} : H^{1}(J_{\ok},\bZ_{l})\simto
H^{1}(C_{\ok},\bZ_{l}),
\\
&f_{*} : H^{1}(C_{\ok},\bZ_{l})\simto
H^{2g-1}(J_{\ok},\bZ_{l}(g-1)).
\endaligned
$$
These are independent of the choice of the
$ k $-rational point of
$ C $,
because a translation on
$ J $ acts trivially on the cohomology of
$ J $.
Since
$ f_{*} $ and
$ f^{*} $ are dual of each other, the canonical pairing on
$ H^{1}(C_{\ok},\bZ_{l}) $ is identified by
$ f^{*} $ with the pairing on
$ H^{1}(J_{\ok},\bZ_{l}) $ given by Poincar\'e duality and
$ f_{*}\scirc f^{*} $.

Let
$ \Gamma_{1} = m^{*}\Theta - pr_{1}^{*}\Theta - pr_{2}^{*}\Theta $,
where
$ \Theta $ is the theta divisor, and
$ m : J \mtim_{k} J \to J $ is the multiplication.
Let
$ \Gamma_{2} \in \CH^{g}(J\mtim_{k}J) $ be the diagonal of
$ f(C) \subset J $ so that
$ f_{*}\scirc f^{*} $ is identified with
$ \Gamma_{2} $.
Since the canonical pairing (3.1.3) can be identified with
Poincar\'e duality (see (3.6)), the assertion is reduced by (3.3)
(where
$ A = X = J $ and
$ D = \Gamma_{1} $) to that the actions of the correspondences
$$
\aligned
&(\Gamma_{1})_{*} :
H^{2g-1}(J_{\ok},\bZ_{l}) \to
H^{1}(J_{\ok},\bZ_{l}(1-g)),
\\
&(\Gamma_{2})_{*} :
H^{1}(J_{\ok},\bZ_{l}(1-g)) \to
H^{2g-1}(J_{\ok},\bZ_{l})
\endaligned
$$
are inverse of each other up to sign, or equivalently, that
the action of the composition of
$ \Gamma_{2}\scirc \Gamma_{1} $ on
$ H^{2g-1}(J_{\ok},\bZ_{l}) $ is the multiplication by
$ -1 $.
Here it is enough to show the assertion for the action on the
Albanese variety of
$ J $, see (3.7).
For
$ a, b \in J(\ok) $, we see that the image of
$ [a] - [b] $ by the action of
$ \Gamma_{2}\scirc\Gamma_{1} $ is given by
$$
f_{*}f^{*}(T_{a}^{*}\Theta - T_{b}^{*}\Theta).
\leqno(3.4.1)
$$

Let
$ C^{(j)} $ denote the
$ j $-th symmetric power of
$ C $.
Then
$ f $ induces
$ f^{(j)} : C^{(j)} \to J $, and
$ f^{(g)} $ is birational [31].
So there is a nonempty Zariski-open subset
$ U $ of
$ J $ such that for
$ a \in U(\ok) $,
there exists uniquely
$ \{c_{1}, \dots , c_{g}\} \in C^{(g)}(\ok) $ satisfying
$$
-a -\sum_{i\ne j}f(c_{i}) = f(c_{j}).
\leqno(3.4.2)
$$
Since
$ \Theta = f^{(g-1)}(C^{(g-1)}) $,
it implies that
$ {T}_{a}^{*}\Theta^{-} \cap f(C) = \{f(c_{j})\} $,
where
$ \Theta^{-} = (-1)^{*}\Theta $.
But
$ \varphi_{\Theta} = \varphi_{\Theta^{-}} $,
because the action of
$ -1 : J \to J $ on
$ \NS(J) $ is the identity.
So the assertion follows from (3.4.1--2).

\medskip\noindent
{\bf 3.5.~First cohomology case.} Let
$ Y $ be a smooth projective variety with an ample divisor class
$ L $.
Then the canonical pairing on
$ H^{1}(Y_{\ok},\bQ_{l}) $,
defined by Poincar\'e duality and
$ L $,
is abelian-positive.
Indeed, we may assume
$ L $ is very ample, and take
$ C $ a smooth closed subvariety of dimension
$ 1 $ which is an intersection of general hyperplane sections
(replacing
$ k $ with a finite extension if necessary).
Then the composition of the restriction and Gysin morphisms
$$
H^{1}(Y_{\ok},\bQ_{l})\to H^{1}(C_{\ok},\bQ_{l})\to
H^{2n-1}(Y_{\ok},\bQ_{l}(n-1))
$$
coincides with
$ L^{n-1} $,
where
$ n = \dim X $.
So the pairing on
$ H^{1}(Y_{\ok},\bQ_{l}) $ is identified with the restriction of
the natural pairing on
$ H^{1}(C_{\ok},\bQ_{l}) $ to
$ H^{1}(Y_{\ok},\bQ_{l}) $.

\medskip\noindent
{\bf 3.6.~Complement to the proof of (3.4), I.}
Let
$ X $ be an irreducible smooth projective variety over a field
$ k $ having a
$ k $-rational point
$ 0 $.
Let
$ P_{X} $ be the Picard variety of
$ X $.
Then we have a canonical morphism
$ \Alb : X \to P_{X}^{\vee} $ sending
$ 0 $ to
$ 0 $ so that
$ P_{X}^{\vee} $ is identified with the Albanese variety
$ \Alb_{X} $ of
$ X $.
(Indeed, for an abelian variety
$ A $ and a morphism
$ f : (X,0) \to (A,0) $,
we have
$ f^{\vee} : A^{\vee} \to P_{X} $ and
$ f^{\vee \vee} : P_{X}^{\vee} \to A $ so that
$ f^{\vee \vee}\scirc\Alb = f $,
using the theory of divisorial correspondences.)

Let
$ n = \dim X $,
and
$ l $ a prime number different from the characteristic of
$ k $.
Let
$ V_{l}M = T_{l}M\otimes_{\bZ_{l}}\bQ_{l} $ for an abelian
group
$ M $.
Then, using the Kummer sequence together with the arguments in
(3.1), we have canonical isomorphisms
$$
\aligned
V_{l}P_{X}(\ok)
&= H^{1}(X_{\ok},\bQ_{l}(1)),
\\
V_{l}\Alb_{X}(\ok)
&= \Hom(V_{l}P_{X}(\ok),\bQ_{l}(1))
= H^{2n-1}(X_{\ok},\bQ_{l})(n),
\endaligned
$$
where the last isomorphism follows from the first together with
Poincar\'e duality, and the last term is identified with
$ H_{1}(X_{\ok},\bQ_{l})
:= \Hom(H^{1}(X_{\ok},\bQ_{l}),\bQ_{l}) $.

\medskip\noindent
{\bf 3.7.~Complement to the proof of (3.4), II.}
Let
$ X, Y $ be irreducible smooth projective varieties over a field
$ k $ having a
$ k $-rational point
$ 0 $.
Let
$ \Gamma \in \CH^{m}(X\mtim_{k}Y) $ with
$ m = \dim Y $.
Then
$ \Gamma $ induces morphisms
$$
(X,0) \to (\Alb_{Y},0),\quad \Gamma_{*} : \Alb_{X} \to \Alb_{Y},
$$
and the induced morphism
$ \Gamma_{*} : V_{l}\Alb_{X}(\ok) \to V_{l}\Alb_{Y}(\ok) $
is identified by (3.6) with
$$
\Gamma_{*} : H^{2n-1}(X_{\ok},\bQ_{l})(n)
\to H^{2m-1}(Y_{\ok},\bQ_{l})(m),
$$
which is given by
$ cl(\Gamma)^{1,2m-1} \in H^{1}(X_{\ok},\bQ_{l})
\otimes H^{2m-1}(Y_{\ok},\bQ_{l})(m) $, the
$ (1,2m-1) $-component of the cycle class of
$ \Gamma_{*} $, see also [22].
Although this is well-known to specialists, its proof does
not seem to be completely trivial, and we give here a short
sketch for the convenience of the reader.

First we may replace
$ Y $ with
$ \Alb_{Y} $ (by composing
$ \Gamma $ with the graph of
$ Y \to \Alb_{Y}) $ so that the assertion is reduced to the case
$ Y $ is an abelian variety
$ A $.
Consider
$ \Alb(\Gamma) \in \CH^{m}(X\mtim_{k}A) $ which is a graph of a
morphism of
$ X $ to
$ A $,
and is defined by using the additive structure of
$ A $ (applied to the restriction of
$ \Gamma $ to the generic fiber of
$ X\mtim_{k}A \to X) $.
Then the assertion is easily verified if
$ \Gamma $ is replaced by
$ \Alb(\Gamma) $,
because
$ \Alb(\Gamma) $ is then extended to an element of
$ \CH^{m}(\Alb_{X}\mtim_{k}A) $ which is a graph of a morphism of
abelian varieties.

So the assertion is reduced to that
$ cl(\Gamma)^{1,2m-1} = 0 $ if
$ \Alb(\Gamma) = 0 $.
Here we may assume
$ k $ is algebraically closed.
Replacing
$ X $ with a variety which is \'etale over
$ X $,
we may assume that
$ \Gamma $ is a linear combination of the graphs of morphisms of
$ X $ to
$ A $ (where
$ X $ is smooth and irreducible, but may be nonproper).
By induction on the number of the components of
$ \Gamma $,
we may assume that
$$
\Gamma = \Gamma_{g_{1}+g_{2}} - \Gamma_{g_{1}}
- \Gamma_{g_{2}} - \Gamma_{0}
$$
for morphisms
$ g_{i} : X \to A \,(i = 1,2) $,
where
$ \Gamma_{g_{i}} $ denotes the graph of
$ g_{i} $.

Let
$ \tGam $ denote the pullback of the cycle
$ \Gamma_{g_{1}} - \Gamma_{0} $ by the projection
$ A\mtim_{k}X\mtim_{k}A \to X\mtim_{k}A $ sending
$ (a,x,b) $ to
$ (x,a+b) $.
Then
$ \Gamma_{g_{1}} - \Gamma_{0} $ and
$ \Gamma_{g_{1}+g_{2}} - \Gamma_{g_{2}} $ coincide with the
pull-backs of
$ \tGam $ by the inclusions
$ X\mtim_{k}A \to A\mtim_{k}X\mtim_{k}A $ sending
$ (x,a) $ to
$ (0,x,a) $ and
$ (-g_{2}(x),x,a) $ respectively.
Since these inclusions are sections of the projection to the second
and third factors, it is enough to show that the K\"unneth component
of the cycle class of
$ \tGam $ in
$$
H^{1}(A\mtim_{k}X,\bQ_{l})\otimes H^{2m-1}(A,\bQ_{l})(m)
$$
comes from
$ H^{1}(X,\bQ_{l})\otimes H^{2m-1}(A,\bQ_{l})(m) $ by
$ pr_{2}\mtim id $, where
$ pr_{2} : A\mtim_{k}X \to X $ is the second projection.
The last assertion is easily verified by using the K\"unneth
decomposition.

\newpage
\centerline{{\bf 4. Standard Conjectures}}

\bigskip\noindent
{\bf 4.1.}
Let
$ Y $ be an equidimensional smooth projective variety over a field
$ k $.
We fix an ample divisor class
$ L $ of
$ Y $.
Then
$ L $ acts on \'etale cohomology.
Let
$ n = \dim Y $.
By the hard Lefschetz theorem [4], we have
$$
L^{j} : H^{n-j}(Y_{\ok},\bQ_{l})\simto
H^{n+j}(Y_{\ok},\bQ_{l}(j))\,\,\,(j > 0),
$$
which implies the Lefschetz decomposition
$$
H^{j}(Y_{\ok},\bQ_{l})=\mopls_{i\ge 0}\,
L^{i}H^{j-2i}(Y_{\ok},\bQ_{l}(-i))^{\prim}.
$$
This induces a morphism
$$
\Lambda : H^{j}(Y_{\ok},\bQ_{l})\to H^{j-2}(Y_{\ok},\bQ_{l}(-1)),
$$
such that for
$ m \in H^{j}(Y_{\ok},\bQ_{l})^{\prim} $,
we have
$ \Lambda (L^{i}m) = L^{i-1}m $ if
$ i > 0 $,
and
$ 0 $ otherwise.
The standard conjecture
$ B(Y) $ asserts that
$ \Lambda $ is algebraically defined as an action of a correspondence.

We define
$$
A^{j}(Y) = \Im(cl : \CH^{j}(Y)_{\bQ} \to
H^{2j}(Y_{\ok},\bQ_{l}(j))),
\leqno(4.1.1)
$$
so that we have the injective morphism
$ A^{j}(Y) \to H^{2j}(Y_{\ok},\bQ_{l}(j)) $.
Note that
$ A^{j}(Y) $ may apparently depend on the choice of
$ l $.
The standard conjecture
$ A(Y,L) $ asserts the isomorphism
$$
L^{n-2i} : A^{i}(Y) \simto A^{n-i}(Y)\,\,\,(0 \le i < n/2).
$$
This is independent of
$ L $ if the
$ A^{i}(Y) $ are finite dimensional, because it is equivalent to
the equality of the dimensions by the hard Lefschetz theorem for
the \'etale cohomology groups.
The conjecture
$ A(Y,L) $ follows from
$ B(Y) $,
and implies that the Lefschetz decomposition is compatible with the
subspace
$ A^{j}(Y) $ so that
$$
A^{j}(Y) = \mopls_{i\ge 0}\, L^{i}A^{j-i}(Y)^{\prim}.
$$
The standard conjecture of Hodge index type
$ I(Y,L) $ asserts that the pairing
$$
(-1)^{j}\langle L^{n-2j}a,b\rangle \quad \text{for}\,\,
a, b \in A^{j}(Y)^{\prim}
$$
is positive definite for
$ 0 \le j \le n/2 $.

We will denote by
$ D(Y) $ the conjecture which asserts the coincidence of the
homological equivalence and the numerical equivalences for the
cycles on
$ Y $ (i.e. the canonical pairing between
$ A^{j}(Y) $ and
$ A^{n-j}(Y) $ is nondegenerate for
$ A^{j}(Y) $).
This is equivalent to
$ A(Y,L) $ under the assumption
$ I(Y,L) $, and implies the injectivity of
$$
A^{j}(Y)\otimes_{\bQ}\bQ_{l} \to
H^{2j}(Y_{\ok},\bQ_{l}(j)),
\leqno(4.1.2)
$$
and also the independence of
$ A^{j}(Y) $ of the choice of
$ l $.
It is known that
$ D(Y) $ is true for divisors (by Matsusaka), and
$ I(Y,L) $ is true for surfaces (by Segre [28]),
see [14], [16] (and also [10], [19]).

By the Lefschetz decomposition, we have an isomorphism
$$
{}^{*} : H^{n+j}(Y_{\ok},\bQ_{l})\simto
H^{n-j}(Y_{\ok},\bQ_{l}(-j))
$$
such that for
$ m \in H^{i}(Y_{\ok},\bQ_{l}(-a))^{\prim} $,
we have
$$
{}^{*}(L^{a}m) = (-1)^{i(i+1)/2}L^{n-i-a}m.
$$
Combined with Poincar\'e duality, this induces a pairing on
$ H^{\ssbull}(Y_{\ok},\bQ_{l}) $ defined by
$ \langle m,{}^{*}n \rangle $ for
$ m, n \in H^{j}(Y_{\ok},\bQ_{l}) $.

For a nonzero correspondence
$ \Gamma \in A^{n}(Y\mtim_{k}Y) \subset
\End(H^{\ssbull}(Y_{\ok},\bQ_{l})) $,
we define
$ \Gamma' $ to be the composition of
$ {}^{*} $,
$ {}^{t}\Gamma $ and
$ {}^{*} $, where
$ {}^{t}\Gamma $ is the transpose of
$ \Gamma $.
Then
$ \Gamma' $ is algebraic and
$$
\Tr(\Gamma'\scirc \Gamma) > 0,
\leqno(4.1.3)
$$
if
$ B(Y) $ and
$ I(Y\mtim_{k}Y,L\otimes 1+1\otimes L) $ are satisfied,
see [14].
(This is an analogue of the positivity of the Losati involution.)
Note that
$ H^{\ssbull}(Y_{\ok},\bQ_{l})\otimes H^{\ssbull}(Y_{\ok},\bQ_{l}) $
is the direct sum of
$$
S_{i,j} := \sum_{a=0}^{i}\sum_{b=0}^{j}
L^{a}H^{n-i}(Y_{\ok},\bQ_{l})^{\prim}\otimes
L^{b}H^{n-j}(Y_{\ok},\bQ_{l})^{\prim}.
$$
For
$ i = j $, the primitive part of
$ S_{i,i} $ of degree
$ 2n $ for the action of
$ L\otimes 1 + 1\otimes L $ is isomorphic to
$ H^{n-i}(Y_{\ok},\bQ_{l})^{\prim}\otimes
H^{n-i}(Y_{\ok},\bQ_{l})^{\prim} $, and this isomorphism is
compatible with the canonical pairing up to a nonzero
multiple constant.

If
$ B(Y) $ and
$ I(Y\mtim_{k}Y,L\otimes 1+1\otimes L) $ are true, then
$ A^{n}(Y\mtim_{k}Y) $ is a semisimple algebra and there are
projectors to primitive parts (and
$ D(Y\mtim_{k}Y) $ holds), see [14], [15], [18].
We have the projector to each cohomology group by the algebraicity
of the K\"unneth components using the Frobenius morphism [13]
(because
$ k $ is assumed to be a finite field).
By Jannsen [12] it is actually enough to assume
$ D(Y\mtim_{k}Y) $ for the semisimplicity of
$ A^{n}(Y\mtim_{k}Y) $.

We will denote by
$$
\iota : A^{n}(Y\mtim_{k}Y) \to
\End(H^{\ssbull}(Y_{\ok},\bQ_{l}))
$$
the canonical injection induced by (4.1.1).

\medskip\noindent
{\bf 4.2.~Remarks.} (i) Assume
$ D(Y\mtim_{k}Y) $ holds.
Then
$ A := A^{n}(Y\mtim_{k}Y) $ is a direct product of full matrix
algebras over skew fields by [12] (together with Wedderburn's
theorem).
For an element
$ f $ of
$ A $,
there exists an idempotent
$ \pi $ such that
$ \Im\,\iota(\pi) = \Im\,\iota(f) $ in
$ H^{\ssbull}(Y_{\ok},{\Bbb Q}_{l}) $.
Indeed, using projectors (and replacing
$ H^{\ssbull}(Y_{\ok},{\Bbb Q}_{l}) $ with the corresponding
subspace
$ V) $, the assertion is reduced to the case where
$ A $ is a full matrix algebra over a skew field.
Then a matrix can be modified to a simple form by using actions of
invertible matrices from both sides as well-known (and a conjugate
of an idempotent is an idempotent).

\medskip
(ii) For a morphism of
$ k $-varieties
$ f : X \to Y $, the idempotent corresponding to the image of
$ f^{*} $ is obtained by considering the graph of
$ f $ and applying the above argument to the disjoint union of
$ X $ and
$ Y $ if
$ \dim X = \dim Y $.
In general we can replace
$ X $ or
$ Y $ by a product with
$ {\bP}^{m} $ (using the K\"unneth decomposition for
$ {\bP}^{m} $), see [27].
Note that
$ D(X\mtim_{k}Y\mtim_{k}{\bP}^{m}) $ can be reduced to
$ D(X\mtim_{k}Y) $.

\medskip\noindent
{\bf 4.3.~Proof of (0.4).}
By (4.1) and (4.2) there exists an idempotent
$ \pi $ of
$ A^{n}(Y\mtim_{k}Y) $ such that
$ \Im\,\iota(\pi) $ coincides with the intersection of the
primitive part with the image of the Cech restriction or co-Cech
Gysin morphism.

We show in general that for a projector
$ \pi $ of
$ A^{n}(Y\mtim_{k}Y) $,
the restriction of the pairing to
$ \Im\,\iota (\pi) $ is nondegenerate if
$ \Im\,\iota (\pi) $ is contained in
$ H^{j}(Y_{\ok},\bQ_{l})^{\prim} $.
For this it is enough to show
$$
\Im\,\iota (\pi) \cap \Ker\,\iota (\pi') = 0,
\leqno(4.3.1)
$$
where
$ \pi' $ is the composition of
$ {}^{*}$,
$ {}^{t}\pi $,
$ {}^{*}$ as in (4.1.3).
We see that (4.3.1) follows from (4.1.3) if
$ \pi $ corresponds to a simple motive, because the intersection is
defined as a motive (and the forgetful functor associating the
underlying vector space commutes with
$ \Im $,
$ \Ker $ and the intersection).
In general, we consider a simple submotive of
$ \Im\,\pi $.
Let
$ \pi_{0} $ be the projector defining it.
Then (4.3.1) holds for
$ \pi_{0} $,
and we get a decomposition
$$
H^{j}(Y_{\ok},\bQ_{l})^{\prim} = \Im\,\iota (\pi_{0}) \oplus
(\Ker\,\iota (\pi'_{0}) \cap H^{j}(Y_{\ok},\bQ_{l})^{\prim}),
$$
which is defined motivically, and is compatible with
$ \Im\,\iota (\pi) $, i.e.
$$
\Im\,\iota (\pi) = \Im\,\iota (\pi_{0}) \oplus
(\Ker\,\iota (\pi'_{0}) \cap \Im\,\iota (\pi)).
$$
Therefore, replacing
$ H^{j}(Y_{\ok},\bQ_{l})^{\prim} $ with
$ \Ker\,\iota (\pi'_{0}) \cap H^{j}(Y_{\ok},\bQ_{l})^{\prim} $,
we can proceed by induction.
This completes the proof of (0.4).

\bigskip\bigskip\centerline{\bf 5. Frobenius Action}

\bigskip\noindent
{\bf 5.1.~Weil Conjecture.}
Let
$ Y $ be an equidimensional smooth projective
$ k $-variety, and put
$ A = A^{n}(Y\mtim_{k}Y) $ where
$ n = \dim Y $.
Let
$ q = |k| $,
$ p = \char k $.
We denote by
$ g $ the graph of the
$ q $-th power Frobenius.
Then
$ g $ belongs to the center of
$ A $ (using the compatibility with the Galois action).
To simplify the relation with the eigenvalues, let
$ P(T) $ denote the characteristic polynomial of the action of
the Frobenius
$ g $ on
$ H^{\ssbull}(Y_{\ok},\bQ_{l}) \,(:= \mopls_{i\in \bZ}
H^{i}(Y_{\ok},\bQ_{l})) $ such that the eigenvalues of the
Frobenius action are the roots of
$ P(T) $ (this normalization is different from the one used in
[4], [13], etc.)
Then
$ P(T) $ is a monic polynomial with integral coefficients and is
independent of
$ l \,(\ne p) $.
The eigenvalues of the action of
$ g $ on
$ H^{i}(Y_{\ok},\bQ_{l}) $ are algebraic integers whose image
by any embedding of
$ \obQ $ into
$ \bC $ has absolute value
$ q^{i/2} $,
see [4].

Let
$ E $ be an algebraic number field.
Consider the decomposition
$ P(T) = \prod_{j} P_{j}(T)^{m_{j}} $ in
$ E[T] $ where the
$ P_{j}(T) $ are monic irreducible polynomials whose
coefficients are algebraic integers in
$ E $.
Let
$ v $ be a prime (i.e. a finite place) of
$ E $ over a given prime number
$ l \,(\ne p) $.
Let
$ E_{v} $ denote the completion of
$ E $ at
$ v $,
and
$ \oE_{v} $ denote an algebraic closure of
$ E_{v} $ (which is isomorphic to
$ \obQ_{l}) $.
Consider the ring of correspondences
$ A_{E} $ with
$ E $-coefficients.
We have a natural inclusion
$$
\iota : A_{E} \to \End(H^{\ssbull}(Y_{\ok},E_{v})),
$$
and a natural surjection
$ A\otimes_{\bQ}E \to A_{E} $,
but the injectivity of the last morphism is not clear unless
the conjecture
$ D(Y\mtim_{k}Y) $ holds.

As well-known, there exist
$ R_{j}(T) \in E[T] $ such that
$ \pi_{j} := R_{j}(g) \in A_{E} $ is an idempotent,
$ \sum_{j} \pi_{j} = 1 $ and the characteristic polynomial
of the action of
$ g $ on
$ \Im\, \iota (\pi_{j}) $ is
$ P_{j}(T)^{m_{j}} $.
Note that
$ \Im\, \iota (\pi_{j}) $ is contained in some cohomology group
$ H^{i}(Y_{\ok},E_{v}) $ by [4], see [13].

\medskip\noindent
{\bf 5.2.~Proposition.} {\it
Let
$ M_{j} $ be the motive defined by
$ \pi_{j} $.
Assume
$ m_{j} = 1 $.
Then the endomorphism ring
$ \End(M_{j}) $ is generated by
$ g $ over
$ E $, and is isomorphic to
$ E[T]/(P_{j}(T)) $.
}

\medskip\noindent
{\it Proof.}
Let
$ F $ be a minimal Galois extension of
$ E $ containing all the roots of
$ P_{j}[T] $.
Let
$ d_{j} = \deg P_{j} $.
For a root
$ \alpha $ of
$ P_{j}[T] $ in
$ F $,
let
$ E[\alpha ] $ denote the subfield of
$ F $ generated by
$ \alpha $, which is isomorphic to
$ E[T]/(P_{j}(T)) $.
By the same argument as above, we have a projector
$ \pi_{\alpha} $ in
$ A_{E[\alpha]} $ such that
$$
g\pi_{\alpha} = \alpha \pi_{\alpha}.
\leqno(5.2.1)
$$
Take any
$ \Gamma \in A_{E} $.
There exists a polynomial
$ h(T) \in E[T] $ such that
$$
\Gamma \pi_{\alpha}\cdot \Delta = h(\alpha) \in E[\alpha],
\leqno(5.2.2)
$$
where
$ \Delta $ is the diagonal cycle of
$ Y $,
and the left-hand side denotes the intersection number.
Note that
$ h $ is independent of
$ \alpha $ using the action of
$ \Gal(F/E) $ because, fixing one root
$ \beta $,
we have for each root
$ \alpha $ an automorphism
$ \sigma $ of
$ F/E $ such that
$ \alpha = \beta^{\sigma} $.
By the Lefschetz trace formula,
$ h(\alpha) $ coincides with the eigenvalue of the action of
$ \Gamma $ on the
$ \alpha $-eigenspace of
$ g $ up to a sign independent of
$ \alpha $.
Here
$ \alpha $ denotes also the image of
$ \alpha $ by the embedding
$ F \to F_{v} $ where
$ v $ is a prime of
$ F $ over
$ l $.
So the assertion follows.

\medskip\noindent
{\bf 5.3.~Image of correspondences.}
Let
$ X $ be an equidimensional smooth projective
$ k $-variety, and
$ \Gamma \in \CH^{r}(X\mtim_{k}Y)_{E} $.
Assume that
$ E $ is a subfield of
$ \bR $,
$ \deg P_{j} \le 2 $ for any
$ j $, and the roots of
$ P_{j} $ are not real if
$ \deg P_{j} = 2 $.
Let
$ \beta $ be an eigenvalue in
$ \oE_{v} $ of the Frobenius action on
$ \Im\, \Gamma_{*} \cap H^{i}(Y_{\ok},E_{v})^{\prim} $,
and assume that its multiplicity as an eigenvalue of the
Frobenius action on
$ H^{i}(Y_{\ok},E_{v})^{\prim} $ is
$ 1 $ (i.e. the dimension of the generalized eigenspace in
$ H^{i}(Y_{\ok},E_{v})^{\prim} $ is
$ 1 $).
By definition
$ \beta $ is a root of some
$ P_{j}(T) $.
Assume further that
$ P_{j}(T) $ is irreducible over
$ E_{v} $,
and
$ \deg P_{j} = 2 $ (because the case
$ \deg P_{j} = 1 $ is trivial).
Let
$ \beta' $ be the other root of
$ P_{j}(T) $.
Then
$ \beta' $ is an eigenvalue of the Frobenius action on
$ \Im\, \Gamma_{*} \cap H^{i}(Y_{\ok},E_{v})^{\prim} $,
and hence
$ \Im\, \Gamma_{*} \cap H^{i}(Y_{\ok},E_{v})^{\prim} $ contains
$ \Im\, \iota (\pi_{j}) $ on which the canonical pairing is
nondegenerate.
Indeed, the pairing is compatible with the
Frobenius action and has value in
$ E_{v}(-i) $ on which the action of
$ g $ is a multiplication by
$ q^{i} $.
So we are interested in the problem: When is
$ P_{j} $ irreducible over
$ E_{v} $?

Let
$ F $ be the field generated by all the roots
$ \alpha $ of
$ P $ in
$ \bC $.
Since
$ \alpha \oal = q^{r} $ with
$ r \in \bZ $,
the complex conjugation on
$ F $ (induced by that on
$ \bC/\bR) $ is in the center of
$ \Gal(F/\bQ) $,
because
$ \alpha^{\sigma}\oal^{\sigma} = q^{r} $
and hence
$ \oal^{\sigma} = \overline{\alpha^{\sigma}} $ for
$ \sigma \in \Gal(F/\bQ) $.
Let
$ E $ be the subfield of
$ F $ fixed by the complex conjugation.
By the Tchebotarev density theorem, there are infinitely many
primes
$ v $ whose density is
$ 1/2 $ and such that
$ P_{j} $ modulo
$ v $ is irreducible over the residue field (hence
$ P_{j} $ is irreducible over
$ E_{v} $) considering the conjugacy class of the above
complex conjugation in
$ \Gal(F/E) $.
So we get the following.

\medskip\noindent
{\bf 5.4.~Proposition.} {\it
For prime numbers
$ l $ such that some
$ v $ as above is over
$ l $, Conjecture {\rm (0.1)} is true if each eigenvalue of the
Frobenius action on the intersections of the primitive part with
$ \Im\, \rho $ and with
$ \Im\, \gamma $ has multiplicity
$ 1 $ as an eigenvalue of the Frobenius action on the primitive
cohomology
$ H^{j}(Y^{(i+1)},\bQ_{l})^{\prim} $ for
$ j \ge 2 $ and
$ i \ge 0 $.
}

\medskip\noindent
{\bf 5.5.~Remark.}
If Conjecture (0.1) holds for a general hyperplane section of the
generic fiber, then it is enough to consider the intersection with
$ \Im\, \rho $ by Theorem (0.3).
For a general
$ l $, Conjecture (0.1) in the case of (5.4) can be reduced to the
conjecture
$ D $ (or
$ B $).
Thus we can avoid the positivity in this simple case.
However, it is not easy to avoid it even in the case where the
multiplicity
$ 1 $ holds for each irreducible component of
$ Y^{(i+1)} $.
Indeed, the ambiguity of the pairing on the motive
$ M_{j} $ in (5.2) is given by an automorphism
$ h \in \bQ[g] \,(= \bQ[T]/(P_{j}(T))) $ where
$ E = \bQ $.
If there is a morphism of
$ M_{j} $ to a direct sum of simple motives which are isomorphic to
$ M_{j} $ and indexed by
$ r $, and if the morphism to the
$ r $-th factor is given by a correspondence
$ \Gamma_{r} $, then the pull-back of the pairing corresponds to the
sum of
$ h_{r} := \Gamma'_{r} \scirc \Gamma_{r} \in \bQ[g] $ up to a sign,
and the problem is closely related to the standard conjecture of
Hodge index type.

\bigskip\bigskip
\centerline{{\bf References}}

\bigskip

\item{[1]}
A. Beilinson, J. Bernstein and P. Deligne, Faisceaux pervers,
Ast\'erisque, vol. 100, Soc. Math. France, Paris, 1982.

\item{[2]}
A.J. de Jong, Smoothness, semi-stability and alternations, Publ. Math.
IHES, 83 (1996), 51--93.

\item{[3]}
P. Deligne, Th\'eorie de Hodge I, Actes Congr\`es Intern. Math.,
1970, vol. 1, 425-430; II, Publ. Math. IHES, 40 (1971), 5--57;
III ibid., 44 (1974), 5--77

\item{[4]}
\SameAuthor, La conjecture de Weil II, Publ. Math. IHES, 52 (1980),
137--252.

\item{[5]}
\SameAuthor, Positivit\'e: signe II (manuscrit, 6--11--85).

\item{[6]}
\SameAuthor, Dualit\'e, in SGA 4 1/2, Lect. Notes in Math., vol. 569
Springer, Berlin, 1977, pp. 154--167.

\item{[7]}
P. Griffiths and J. Harris, On the Noether-Lefschetz theorem and
some remarks on codimension-two cycles, Math. Ann. 271 (1985),
31--51.

\item{[8]}
A. Grothendieck et al, Groupes de monodromie en g\'eom\'etrie
alg\'ebrique, SGA 7 I, Lect. Notes in Math. vol. 288 Springer, Berlin,
1972.

\item{[9]}
F. Guill\'en and V. Navarro Aznar, Sur le th\'eor\`eme
local des cycles invariants, Duke Math. J. 61 (1990), 133--155.

\item{[10]}
R. Hartshorne, Algebraic Geometry, Springer, New York, 1977.

\item{[11]}
L. Illusie, Autour du th\'eor\`eme de monodromie locale,
Ast\'erisque 223 (1994), 9--57.

\item{[12]}
U. Jannsen, Motives, numerical equivalence, and semi-simplicity, Inv.
Math. 107 (1992), 447--452,

\item{[13]}
N. Katz and W. Messing, Some consequences of the Riemann hypothesis
for varieties over finite fields, Inv. Math. 23 (1974), 73--77.

\item{[14]}
S. Kleiman, Algebraic cycles and Weil conjecture, in Dix expos\'es
sur la cohomologie des sch\'emas, North-Holland, Amsterdam, 1968, pp.
359--386.

\item{[15]}
\SameAuthor, Motives, in Algebraic Geometry (Oslo, 1970),
Wolters-Noordhoff, Groningen,
1972, pp. 53--82.

\item{[16]}
\SameAuthor, The standard conjecture, Proc. Symp. Pure Math. 55
(1994), Part 1, 3--20.

\item{[17]}
S. Lang, Abelian varieties, Interscience Publishers, New York, 1959.

\item{[18]}
Y.I. Manin, Correspondences, motifs and monoidal transformations,
Math. USSR Sb. 6 (1968), 439--470.

\item{[19]}
J.S. Milne, Etale cohomology, Princeton University Press, 1980.

\item{[20]}
D. Mumford, Abelian varieties, Oxford University Press, 1970.

\item{[21]}
M. Rapoport and T. Zink,
\"Uber die lokale
Zetafunktion von Shimuravariet\"aten. Monodromiefiltration und
verschwindende Zyklen in ungleicher Charakteristik, Inv. Math. 68
(1982), 21--101.

\item{[22]}
W. Raskind, Higher
$ l $-adic Abel-Jacobi mappings and filtrations on Chow groups,
Duke Math. J. 78 (1995), 33--57.

\item{[23]}
M. Saito, Modules de Hodge polarisables, Publ. RIMS, Kyoto Univ., 24
(1988), 849--995.

\item{[24]}
M. Saito and S. Zucker, The kernel spectral sequence of vanishing
cycles, Duke Math. J. 61 (1990), 329--339.

\item{[25]}
T. Saito, Weight monodromy conjecture for
$ l $-adic representations associated to modular forms, in the
Arithmetic and Geometry of Algebraic Cycles, Kluwer Academic,
Dordrecht, 2000, pp. 427--431.

\item{[26]}
W. Schmid, Variation of Hodge structure: the singularities of the
period mapping, Inv. Math. 22 (1973), 211--319.

\item{[27]}
A. J. Scholl, Classical Motives, Proc. Symp. Pure Math.
55 (1994), Part 1, 163--187.

\item{[28]}
B. Segre, Intorno ad teorema di Hodge sulla teoria della base per
le curve di una superficie algebrica, Ann. Mat. 16 (1937), 157--163.
	
\item{[29]}
J.H.M. Steenbrink, Limits of Hodge structures, Inv. Math. 31
(1975/76), no. 3, 229--257.

\item{[30]}
J. Tate, Conjectures on algebraic cycles in
$ l $-adic cohomology, Proc. Symp. Pure Math. 55 (1994), Part 1,
71--83.

\item{[31]}
A. Weil, Vari\'et\'es Ab\'eliennes et Courbes Alg\'ebriques,
Hermann, Paris, 1948.

\bigskip\noindent
\ver

\bye